\documentclass[a4paper,12pt]{elsarticle}
\usepackage{lineno,hyperref}
\usepackage[latin5]{inputenc}
\usepackage{amsmath}
\usepackage{graphicx}
\usepackage{subfigure}
\usepackage{caption}
\usepackage{float}
\usepackage{geometry}
\begin{document}
\begin{frontmatter}
\title{Exponential B-spline Collocation Solutions to the Gardner Equation}
\author{Ozlem Ersoy Hepson$^{a,*}$, Alper Korkmaz$^b$, Idiris Dag$^c$}
\ead{ozersoy@ogu.edu.tr}
\address{$^a$ Eskisehir Osmangazi University, Mathematics \& Computer Department, Eskisehir, Turkey, \\
$^b$Çankırı Karatekin University, Department of Mathematics, Çankırı, Turkey, \\
$^c$ Eskisehir Osmangazi University, Computer Engineering Department, Eskisehir, Turkey.}
\begin{abstract}
The exponential B-spline basis function set is used to develop a collocation method for some initial boundary value problems (IBVPs) to the Gardner equation. The Gardner equation has two nonlinear terms, namely quadratic and cubic ones. The order reduction of the equation is resulted in a coupled system of PDEs that enables the exponential B-splines to be implemented. The system is integrated in time by Crank-Nicolson implicit method. The validity of the method is investigated by calculating the discrete maximum error norm and observing the absolute relative changes of the conservation laws at the end of the simulations.
\end{abstract}


\end{frontmatter}

\section{Introduction}
The Gardner equation 
\begin{equation}
u_{t}+\mu_1 uu_{x}+\mu_2 u^2u_{x}+\mu_3 u_{xxx}=0 \label{gardner}
\end{equation}
where $\mu_i,i=1,2,3$ are constant parameters has two nonlinear terms besides the third order dissipative term. When these two nonlinear terms are considered, the equation is a combined form of the Korteweg-de Vries and the modified Korteweg-de Vries equations. The Gardner equation is a significant model for the motion of negative ion-acoustic plasma waves \cite{ruderman1}. The equation can also describe internal waves with large amplitudes and weakly nonlinear dispersive waves \cite{kamc1}. Moreover, undular bore type solutions for the Gardner equation and the relation of these solutions with the sign of the cubic nonlinear term are studied in details by Kamchatnov et al. in the same paper. Various patterns covering undular bores in the bright and dark cnoidals, trigonometric bores, kinks and some other types constructed from them are developed in \cite{kamc1}. Occurrence of large ocean waves in unexpected meaning can be modeled by in the Gardner equation and in some cases the modulational instability can be observed \cite{grimshaw1}. An analytical study deals with transcritical flow of a fluid passing a local topographical obstacle \cite{kamc2}.

Interaction of various types of solitons and other wave profiles are examined in details in different studies. Collision of a large amplitude soliton with a limiting soliton is derived by the aid of the Darboux transformation \cite{sul1}. A large class of interaction solutions of solitons with cnoidal and periodic waves is constructed by the consistent $\tanh$-expansion method \cite{hu1}. Some collusion models of cnoidal wave to soliton are derived by the consistent Riccati expansion approach \cite{feng1}. Some exact traveling wave type solutions constructed in terms of $\tanh$ and $\coth$ are constructed by implementing the extended form of the $\tanh$ method \cite{bekir1}. Projective Riccati equations also enable to generate some hyperbolic type solitary wave and periodic solutions expressed in the finite series form \cite{fu1}. These type solutions can also be constructed by generalized form of the $(G'/G)-$expansion \cite{lu1,naher1,taghizade1}, $\exp -$ function \cite{akbar1} methods. A bunch of traveling wave type exact solutions to the Gardner equation are determined by using various hyperbolic ansatzes\cite{wazwaz1}. Zayed and Abdelaziz \cite{zayed1} focus on new exact traveling wave solutions constructed by $(G'/G,1/G)$ approach. Mapping method can also be used to integrate the Gardner equation \cite{krishan1}. Lie group and tan-cot methods are also effective methods to construct some exact solutions to the Gardner equation \cite{guo1,jawad1}.

The three conservation laws representing various physical quantities such as momentum, energy and etc. are  
given for the generalized Gardner equation power law nonlinearities by using some algebraic and derivative manipulations \cite{hamdi1}. Some finite difference and restrictive Taylor's approximation are used to determine the numerical solutions to the Gardner equation \cite{nishiyama1,rageh1}.
In the present study, we construct exponential B-spline collocation algorithms for the numerical solutions to some IBVPs for the Gardner equation. Even though the cubic degree of the basis does not allow to approach the third order derivative, the reduction of the derivative order permits to implement of the method to the resultant nonlinear PDE system.
 
Let $v=u_x$. Then, the Gardner equation (\ref{gardner}) is reduced to a coupled system nonlinear PDEs
\begin{equation}
\begin{aligned}
u_{t}+\left(\mu_1 u+\mu_2 u^{2}\right)u_x+\mu_3 v_{xx}&=0 \\ 
v-u_{x}&=0
\end{aligned}
\label{2}
\end{equation}
The initial data
\begin{equation}
\begin{aligned}
u(x,0)&=\omega(x) \\
v(x,0)&=\omega_x(x)
\end{aligned}
\end{equation}
and the zero Neumann conditions at both end of the artificial finite problem interval $[a,b]$ are combined with the system (\ref{2}) for the mathematical representation of the IBVPs.
\section{Exponential B-spline Approach}
Let $\pi$ be a uniformly distributed grids of the finite interval $[a,b]$, such as,
$$\pi : x_m=a+mh, m=0,1, \ldots N$$
where $h=\frac{b-a}{N}$. Then, the exponential cubic B-splines set  
\begin{equation}
B_{m}(x)=\left\{ 
\begin{array}{lcc}
b_{2}\left ( (x_{m-2}-x)-\frac{1}{\zeta} \sinh(\zeta(x_{m-2}-x)) \right ) & , & [x_{m-2},x_{m-1}] \\ 
a_{1}+b_{1}(x_{m}-x)+c_{1}\exp(\zeta(x_{m}-x))+d_{1}\exp(-\zeta(x_{m}-x)) & , & [x_{m-1},x_{m}] \\ 
a_{1}+b_{1}(x-x_{m})+c_{1}\exp(\zeta(x-x_{m}))+d_{1}\exp(-\zeta(x-x_{m})) & , & [x_{m},x_{m+1}] \\ 
b_{2}\left ( (x-x_{m+2})-\frac{1}{\zeta} \sinh(\zeta(x-x_{m+2})) \right ) & , & [x_{m+1},x_{m+2}] \\ 
0 & , & otherwise%
\end{array}%
\right.\label{ecbs}
\end{equation}%
where $m=-1,0, \ldots N+1$, 
$$
\begin{array}{l}
a_{1}=\dfrac{\zeta h\cosh(\zeta h)}{\zeta h\cosh(\zeta h)-\sinh(\zeta h)}, \\
b_{1}=\dfrac{\zeta }{2} \dfrac{\cosh(\zeta h)(\cosh(\zeta h)-1)+\sinh^2(\zeta h)}{(\zeta h\cosh(\zeta h)-\sinh(\zeta h))(1-\cosh(\zeta h))}, \\
b_{2}=\dfrac{\zeta }{2(\zeta h\cosh(\zeta h)-\sinh(\zeta h))}, \\
c_{1}=\dfrac{1}{4} \dfrac{\exp(-\zeta h)(1-\cosh(\zeta h))+\sinh(\zeta h)(\exp(-\zeta h)-1))}{(\zeta h\cosh(\zeta h)-\sinh(\zeta h))(1-\cosh(\zeta h))}, \\
d_{1}=\dfrac{1}{4} \dfrac{\exp(\zeta h)(\cosh(\zeta h)-1)+\sinh(\zeta h)(\exp(\zeta h)-1))}{(\zeta h\cosh(\zeta h)-\sinh(\zeta h))(1-\cosh(\zeta h))},
\end{array}
$$ 
with real parameter $\zeta$ constitutes a basis functions set for the functions defined over the problem interval $[a,b]$ \cite{mccartin}. Each exponential cubic B-spline $B_m(x)$ has two continuous first and second order derivatives defined in the interval $[x_{m-2},x_{m+2}]$. Exponential B-splines are used as basis functions in some methods to solve problems appearing in various fields \cite{moh1,moh2,ozkdv,ozsivas,ozfish,ozreact}. The functional and derivative values of the exponential B-splines are summarized in Table \ref{table:ebs} .

\scriptsize{
\begin{table}[h]
\caption{$B_{m}(x)$ and its first and second derivatives}
\scriptsize{
\begin{tabular}{lccccc}
\hline \hline
$x$ & $x_{m-2}$ & $x_{m-1}$ & $x_{m}$ & $x_{m+1}$ & $x_{m+2}$ \\ 
\hline \\
$B_{m}$ & $0$ & $\dfrac{\sinh(\zeta h)-\zeta h}{2(\zeta h\cosh(\zeta h)-\sinh(\zeta h))}$ & $1$ & $\dfrac{\sinh(\zeta h)-\zeta h}{2(\zeta h\cosh(\zeta h)-\sinh(\zeta h))}$ & $%
0$ \\ 
$B_{m}^{^{\prime }}$ & $0$ & $\dfrac{\zeta (1-\cosh(\zeta h))}{2(\zeta h\cosh(\zeta h)-\sinh(\zeta h))}$ & $0$ & $\dfrac{%
\zeta (\cosh(\zeta h)-1)}{2(\zeta h\cosh(\zeta h)-\sinh(\zeta h))}$ & $0$ \\ 
$B_{m}^{^{\prime \prime }}$ & $0$ & $\dfrac{\zeta ^{2}\sinh(\zeta h)}{2(\zeta h\cosh(\zeta h)-\sinh(\zeta h))}$ & $-\dfrac{%
\zeta ^{2}\sinh(\zeta h)}{\zeta h\cosh(\zeta h)-\sinh(\zeta h)}$ & $\dfrac{\zeta ^{2}\sinh(\zeta h)}{2(\zeta h\cosh(\zeta h)-\sinh(\zeta h))}$ & $0$ \\ 
\hline \hline
\end{tabular}}
\label{table:ebs}
\end{table}}
\normalsize
Let $U$ and $V$ be approximate solutions to $u$ and $v$ respectively. Then,

\begin{equation}
U(x,t)=\sum_{m=-1}^{N+1}\delta _{m}B_{m}(x),\text{ }V(x,t)=\sum_{m=-1}^{N+1}%
\phi _{m}B_{m}(x).  \label{7}
\end{equation}%
where $\delta _{m}$ and $\phi _{m}$ are time dependent parameters to be
determined from the collocation points $x_{m},m=0,...,N$ and the complementary data.
 The nodal values $U$ and its first and second
derivatives at the knots can be found from the (\ref{7}) as 
\begin{eqnarray}
&&%
\begin{tabular}{l}
$U_{m}=U(x_{m},t)=\dfrac{s-\zeta h}{2(\zeta hc-s)}\delta _{m-1}+\delta _{m}+\dfrac{s-\zeta h%
}{2(\zeta hc-s)}\delta _{m+1},$ \\ 
$U_{m}^{\prime }=U^{\prime }(x_{m},t)=\dfrac{\zeta (1-c)}{2(\zeta hc-s)}\delta _{m-1}+%
\dfrac{\zeta (c-1)}{2(\zeta hc-s)}\delta _{m+1}$ \\ 
$U_{m}^{\prime \prime }=U^{\prime \prime }(x_{m},t)=\dfrac{\zeta ^{2}s}{2(\zeta hc-s)}%
\delta _{m-1}-\dfrac{\zeta ^{2}s}{\zeta hc-s}\delta _{m}+\dfrac{\zeta ^{2}s}{2(\zeta hc-s)}%
\delta _{m+1}$%
\end{tabular}
\label{f5} \\
&&%
\begin{tabular}{l}
$V_{m}=V(x_{m},t)=\dfrac{s-\zeta h}{2(\zeta hc-s)}\phi _{m-1}+\phi _{m}+\dfrac{s-\zeta h}{%
2(\zeta hc-s)}\phi _{m+1},$ \\ 
$V_{m}^{\prime }=V^{\prime }(x_{m},t)=\dfrac{\zeta (1-c)}{2(\zeta hc-s)}\phi _{m-1}+%
\dfrac{\zeta (c-1)}{2(\zeta hc-s)}\phi _{m+1}$ \\ 
$V_{m}^{\prime \prime }=V^{\prime \prime }(x_{m},t)=\dfrac{\zeta ^{2}s}{2(\zeta hc-s)}%
\phi _{m-1}-\dfrac{\zeta ^{2}s}{\zeta hc-s}\phi _{m}+\dfrac{\zeta ^{2}s}{2(\zeta hc-s)}\phi
_{m+1}.$%
\end{tabular}
\label{f6}
\end{eqnarray}%
When Gardner equation is space-splitted as (\ref{2}), the resultant system has
maximum second order derivatives that enables to construct a smooth cubic B-spline approximation 
with the exponential B-splines. Implementation of the Crank-Nicolson method to the space-slitted system (\ref{2}) yields
\begin{equation}
\begin{array}{r}
\dfrac{U^{n+1}-U^{n}}{\Delta t}+\mu _{1}\dfrac{(UU_x)^{n+1}+(UU_x)^{n}}{2}+\mu
_{2}\dfrac{(U^{2}U_x)^{n+1}+(U^{2}U_x)^{n}}{2}+\mu _{3}\dfrac{%
V_{xx}^{n+1}+V_{xx}^{n}}{2}=0 \\ 
\\ 
\dfrac{V^{n+1}+V^{n}}{2}-\dfrac{U_{x}^{n+1}+U_{x}^{n}}{2}=0%
\end{array}
\label{10}
\end{equation}
where $t^{n+1}=t^{n}+\Delta t$, and the superscripts $n$ and $n+1$ denote $(n)$th and $(n+1)$th time levels, respectively.

The nonlinear terms $(UU_x)^{n+1}$ and $\left( U^{2}U_x\right) ^{n+1}$ in Eq. (\ref{10}) are converted to linear forms by using 
\begin{equation*}
(UU_x)^{n+1}=U^{n+1}U_x^{n}+U^{n}U_x^{n+1}-U^{n}U_x^{n}
\end{equation*}%
and%
\begin{equation*}
(U^{2}U_x)^{n+1}=2U^{n+1}U^{n}U_x^{n}+(U^{n})^{2}U_x^{n+1}-2(U^{n})^{2}U_x^{n}
\end{equation*}
defined in \cite{rubin}. The resultant linear system is discretized in time by using Crank-Nicolson method as\scriptsize{
\begin{equation}
\begin{aligned}
&\left[ \left( \frac{2}{\Delta t}+\mu _{1}L+2\mu _{2}KL\right) \alpha _{1}+\left( \mu _{1}K+\mu _{2}K^{2}\right)
\beta _{1}
\right] \delta _{j-1}^{n+1}+\left[ \mu _{3}\gamma _{1}\right] \phi _{j-1}^{n+1}+\left[ \left( \frac{%
2}{\Delta t}+\mu _{1}L+2\mu _{2}KL\right) \alpha _{2}\right] \delta
_{j}^{n+1}\\
&+\left[ \mu_{3}\gamma _{2}\right] \phi _{j}^{n+1} \left[ \left( \frac{2}{\Delta t}+\mu _{1}L+2\mu _{2}KL\right) \alpha _{1}-\left( \mu _{1}K+\mu _{2}K^{2}\right)
\beta _{1}\right] \delta _{j+1}^{n+1}+\left[ \mu _{3}\gamma _{1}\right] \phi _{j+1}^{n+1}\\
&=\left[ \left( \frac{2}{\Delta t}+\mu _{2}KL\right) \alpha _{1}\right]
\delta _{j-1}^{n}-\mu _{3}\gamma _{1}\phi _{j-1}^{n}+\left[ \left( \frac{2}{%
\Delta t}+\mu _{2}KL\right) \alpha _{2}\right] \delta _{j}^{n}-\mu
_{3}\gamma _{2}\phi _{j}^{n}+\left[ \left( \frac{2}{\Delta t}+\mu
_{2}KL\right) \alpha _{1}\right] \delta _{j+1}^{n}&-\mu _{3}\gamma _{1}\phi
_{j+1}^{n} \\
&-\beta _{1}\delta _{j-1}^{n+1}+\alpha _{1}\phi _{j-1}^{n+1}+\alpha _{2}\phi
_{j}^{n+1}+\beta _{1}\delta _{j+1}^{n+1}+\alpha _{1}\phi _{j+1}^{n+1}
=\beta _{1}\delta _{j-1}^{n}-\alpha _{1}\phi _{j-1}^{n}-\alpha _{2}\phi
_{j}^{n}-\beta _{1}\delta _{j+1}^{n}-\alpha _{1}\phi _{j+1}^{n} \\
&m=0,...,N,\quad n=0,1...,
\end{aligned}\label{12}
\end{equation}}
\normalsize
where%
\begin{equation*}
\begin{array}{l}
K=\alpha _{1}\delta _{m-1}^{n}+\delta _{m}^{n}+\alpha _{1}\delta
_{m+1}^{n} \\ 
L=\alpha _{1}\phi _{m-1}^{n}+\phi _{m}^{n}+\alpha _{1}\phi
_{m+1}^{n}%
\end{array}%
\end{equation*}%
\begin{eqnarray*}
\alpha _{1} &=&\dfrac{s-\zeta h}{2(\zeta hc-s)},\text{ }, \\
\gamma _{1} &=&\dfrac{\zeta ^{2}s}{2(\zeta hc-s)},\text{ }\gamma _{2}=-\dfrac{\zeta ^{2}s}{%
\zeta hc-s} \\
\beta _{1} &=&\dfrac{\zeta (1-c)}{2(\zeta hc-s)},\text{ }\beta _{2}=\dfrac{\zeta (c-1)}{2(\zeta hc-s)}
\end{eqnarray*}%
The system (\ref{12}) can be rewritten in the matrix form as
\begin{equation}
\mathbf{Ax}^{n+1}=\mathbf{Bx}^{n}  \label{13}
\end{equation}%
where%
\begin{equation*}
\mathbf{A=}%
\begin{bmatrix}
\nu _{m1} & \nu _{m2} & \nu _{m3} & \nu _{m4} & \nu _{m5} & \nu _{m2} &  & 
&  &  \\ 
-\beta _{1} & \alpha _{1} & 0 & \alpha _{2} & \beta _{1} & \alpha _{1} &  & 
&  &  \\ 
&  & \nu _{m1} & \nu _{m2} & \nu _{m3} & \nu _{m4} & \nu _{m5} & \nu _{m2} & 
&  \\ 
&  & -\beta _{1} & \alpha _{1} & 0 & \alpha _{2} & \beta _{1} & \alpha _{1}
&  &  \\ 
&  &  & \ddots  & \ddots  & \ddots  & \ddots  & \ddots  & \ddots  &  \\ 
&  &  &  & \nu _{m1} & \nu _{m2} & \nu _{m3} & \nu _{m4} & \nu _{m5} & \nu
_{m2} \\ 
&  &  &  & -\beta _{1} & \alpha _{1} & 0 & \alpha _{2} & \beta _{1} & \alpha
_{1}%
\end{bmatrix}%
\end{equation*}

\begin{equation*}
\mathbf{B=}%
\begin{bmatrix}
\nu _{m6} & -\nu _{m2} & \nu _{m7} & -\nu _{m4} & \nu _{m6} & -\nu _{m2} &  & 
&  &  \\ 
\beta _{1} & -\alpha _{1} & 0 & -\alpha _{2} & -\beta _{1} & -\alpha _{1} & 
&  &  &  \\ 
&  & \nu _{m6} & -\nu _{m2} & \nu _{m7} & -\nu _{m4} & \nu _{m6} & -\nu _{m2} & 
&  \\ 
&  & \beta _{1} & -\alpha _{1} & 0 & -\alpha _{2} & -\beta _{1} & -\alpha
_{1} &  &  \\ 
&  &  & \ddots & \ddots & \ddots & \ddots & \ddots & \ddots &  \\ 
&  &  &  & \nu _{m6} & -\nu _{m2} & \nu _{m7} & -\nu _{m4} & \nu _{m6} & -\nu _{m2} \\ 
&  &  &  & \beta _{1} & -\alpha _{1} & 0 & -\alpha _{2} & -\beta _{1} & 
-\alpha _{1}%
\end{bmatrix}%
\end{equation*}%
and%
\begin{equation*}
\begin{array}{ll}
\nu _{m1}=\left( \frac{2}{\Delta t}+\mu _{1}L+2\mu _{2}KL\right) \alpha _{1}+\left( \mu _{1}K+\mu _{2}K^{2}\right)
\beta _{1}
& \nu _{m6}=\left( \frac{2}{\Delta t}+\mu _{2}KL\right) \alpha _{1} \\ 
\nu _{m2}=\mu _{3}\gamma_{1} &  \nu _{m5}=\left( \frac{2}{\Delta t}+\mu _{1}L+2\mu _{2}KL\right) \alpha _{2}-\left( \mu _{1}K+\mu _{2}K^{2}\right)
\beta _{1}\\ 
\nu _{m3}=\left( \frac{2}{\Delta t}+\mu _{1}L+2\mu _{2}KL\right) \alpha _{2}-\left( \mu _{1}K+\mu _{2}K^{2}\right)
\beta _{1}
& \nu _{m7}=\left( \frac{2}{\Delta t}+\mu _{2}KL\right) \alpha _{2} \\ 
\nu _{m4}=\mu _{3}\gamma
_{2} & 
\end{array}%
\end{equation*}
The system (\ref{13}) consists of $2N+2$ linear equations and $2N+6$ unknown
parameters 
\begin{equation*}
\mathbf{x}^{n+1}=(\delta _{-1}^{n+1},\phi _{-1}^{n+1},\delta _{0}^{n+1},\phi
_{0}^{n+1},\ldots ,\delta _{N+1}^{n+1},\phi _{N+1}^{n+1},).
\end{equation*}%
Additional four constraints are required to have the unique solution to this system. These constraints can be determined by the imposition of the boundary data
\begin{equation}
\begin{aligned}
U_{x}(a,t)&=V_{x}(a,t)=0 \\
U_{x}(b,t)&=V_{x}(b,t)=0
\end{aligned}
\end{equation}
to give
\begin{equation*}
\begin{aligned}
\delta _{-1}&=\delta _{1} \\ 
\phi _{-1}&=\phi _{1} \\ 
\delta _{N-1}&=\delta _{N+1} \\ 
\phi _{N-1}&=\phi _{N+1}%
\end{aligned}%
\end{equation*}

When the unknowns $\delta _{-1},\phi _{-1},\delta _{N+1},\phi
_{N+1},$ are eliminated from the system (\ref{13}), the resultant system is a solvable system of $2N+2$ linear equations with $2N+2$ unknowns. The algorithm derived from the classical Thomas
algorithm is used to solve (\ref{13}) in each time level.

In order to start the iteration algorithm, the parameters $\delta _{m}^{0},\phi _{m}^{0},$ $m=-1,\ldots ,N+1$ are required to be determined from 
\begin{equation*}
\begin{array}{l}
U_{x}(a,0)=0=\delta _{-1}^{0}-\delta _{1}^{0}, \\ 
U(x_{m},0)=\alpha _{1}\delta _{m-1}^{0}+\alpha _{2}\delta _{m}^{0}+\alpha
_{1}\delta _{m+1}^{0}=u(x_{m},0),m=1,...,N-1 \\ 
U_{x}(b,0)=0=\delta _{N-1}^{0}-\delta _{N+1}^{0}, \\ 
V_{x}(a,0)=0=\phi _{-1}^{0}-\phi _{1}^{0} \\ 
V(x_{m},0)=\alpha _{1}\phi _{m-1}^{0}+\alpha _{2}\phi _{m}^{0}+\alpha
_{1}\phi _{m+1}^{0}=v(x_{m},0),m=1,...,N-1 \\ 
V_{x}(a,0)=\phi _{N-1}^{0}-\phi _{N+1}^{0}%
\end{array}%
\end{equation*}

\section{Stability Analysis}
The stability of the method is investigated by performing the Von-Neumann analysis where 
\begin{eqnarray}
\delta_{j}^{n} &=&A_1\xi ^{n}\exp (ij\varphi )  \label{k} \\
\phi _{j}^{n} &=&A_2\xi ^{n}\exp (ij\varphi )  \notag
\end{eqnarray}%
\begin{equation*}
\rho=\frac{\xi ^{n+1}}{\xi ^{n}}
\end{equation*}%
Here, $A_1$ and $A_2$ represent the harmonics amplitude. $k$ is the
mode number, $\rho$ is the amplification factor and $\varphi =kh$. The term $U+U^2$ is assumed as locally constant and replaced $\varepsilon$. Substituting \ref{k} into the system 
\begin{eqnarray}
&&a_{1}\delta_{j-1}^{n+1}+a_{2}\delta_{j}^{n+1}+a_{1}\delta_{j+1}^{n+1}+\frac{\lambda
k\varepsilon }{2}(a_{3}\delta_{j-1}^{n+1}-a_{3}\delta_{j+1}^{n+1})+\frac{k\mu_3 }{2}%
(a_{4}\phi _{j-1}^{n+1}+a_{5}\phi _{j}^{n+1}+a_{4}\phi _{j+1}^{n+1})
\label{k1} \\
&=&a_{1}\delta_{j-1}^{n}+a_{2}\delta_{j}^{n}+a_{1}\delta_{j+1}^{n}-\frac{\lambda
k\varepsilon }{2}(a_{3}\delta_{j-1}^{n}-a_{3}\delta_{j+1}^{n})-\frac{k\mu_3 }{2}%
(a_{4}\phi _{j-1}^{n}+a_{5}\phi _{j}^{n}+a_{4}\phi _{j+1}^{n})  \notag
\end{eqnarray}%
\begin{eqnarray}
&&a_{3}\delta_{j-1}^{n+1}-a_{3}\delta_{j+1}^{n+1}-a_{1}\phi _{j-1}^{n+1}-a_{2}\phi
_{j}^{n+1}-a_{1}\phi _{j+1}^{n+1}  \label{k2} \\
&=&-a_{3}\delta_{j-1}^{n}+a_{3}\delta_{j+1}^{n}+a_{1}\phi _{j-1}^{n}+a_{2}\phi
_{j}^{n}+a_{1}\phi _{j+1}^{n}  \notag
\end{eqnarray}%
gives
\begin{eqnarray*}
&&\xi ^{n+1}\left[ A_1\left( 2a_{1}\cos \varphi +a_{2}\right) +\frac{A_2k\mu_3 }{2}%
\left( 2a_{4}\cos \varphi +a_{5}\right) -i\lambda k\varepsilon A_1a_{3}\sin
\varphi \right]  \\
&=&\xi ^{n}\left[ A_1\left( 2a_{1}\cos \varphi +a_{2}\right) -\frac{A_2k\mu_3 }{2}%
\left( 2a_{4}\cos \varphi +a_{5}\right) -i\lambda k\varepsilon A_1a_{3}\sin
\varphi \right] 
\end{eqnarray*}%
\begin{equation*}
\frac{\xi ^{n+1}}{\xi ^{n}}=\frac{\left[ A_1\left( 2a_{1}\cos \varphi
+a_{2}\right) -\frac{A_2k\mu_3 }{2}\left( 2a_{4}\cos \varphi +a_{5}\right)
-i\lambda k\varepsilon A_1a_{3}\sin \varphi \right] }{\left[ A_1\left(
2a_{1}\cos \varphi +a_{2}\right) +\frac{A_2k\mu_3 }{2}\left( 2a_{4}\cos \varphi
+a_{5}\right) -i\lambda k\varepsilon A_1a_{3}\sin \varphi \right] }
\end{equation*}%
\begin{equation}
\rho=\frac{\xi ^{n+1}}{\xi ^{n}}=\frac{X_{1}+iY}{X_{2}-iY}  \label{k3}
\end{equation}%
where%
\begin{eqnarray*}
X_{1} &=&A_1\left( 2a_{1}\cos \varphi +a_{2}\right) -\frac{A_2k\mu_3 }{2}\left(
2a_{4}\cos \varphi +a_{5}\right)  \\
X_{2} &=&A_1\left( 2a_{1}\cos \varphi +a_{2}\right) +\frac{A_2k\mu_3 }{2}\left(
2a_{4}\cos \varphi +a_{5}\right)  \\
Y &=&i\lambda k\varepsilon A_1a_{3}\sin \varphi 
\end{eqnarray*}%
and
\begin{eqnarray*}
&&\xi ^{n+1}\left[ -A_2\left( 2a_{1}\cos \varphi +a_{2}\right) -2iA_1a_{3}\sin
\varphi \right]  \\
&=&\xi ^{n}\left[ A_2\left( 2a_{1}\cos \varphi +a_{2}\right) +2iA_1a_{3}\sin
\varphi \right] 
\end{eqnarray*}%
\begin{equation*}
\frac{\xi ^{n+1}}{\xi ^{n}}=\frac{A_2\left( 2a_{1}\cos \varphi +a_{2}\right)
+2iA_1a_{3}\sin \varphi }{-A_2\left( 2a_{1}\cos \varphi +a_{2}\right)
-2iA_1a_{3}\sin \varphi }
\end{equation*}%
\begin{equation}
\rho =\frac{\xi ^{n+1}}{\xi ^{n}}=\frac{X_{3}+iZ}{X_{4}-iZ}  \label{k4}
\end{equation}%
\begin{eqnarray*}
X_{3} &=&A_2\left( 2a_{1}\cos \varphi +a_{2}\right)  \\
X_{4} &=&-A_2\left( 2a_{1}\cos \varphi +a_{2}\right)  \\
Z &=&2iA_1a_{3}\sin \varphi 
\end{eqnarray*}
It can be concluded from both (\ref{k3}) and (\ref{k4}) that $\left \vert \rho \right \vert$ is less than or equal to $1$. Thus, the proposed method method is unconditionally stable.

\section{Illustrations}

This section contains implementation of the proposed algorithm to some IVPs to validate its accuracy and efficiency. The accuracy of the numerical results is measured by using the discrete maximum norm defined as
\begin{equation*}
\begin{aligned}
L_{\infty}(t)&=\left \vert u(x,t)-U(x,t)\right \vert _{\infty }=\max \limits_{m}\left
\vert u(x_m,t)-U(x_m,t)\right \vert \\
\end{aligned}%
\end{equation*}
at the time $t$ when the analytical solution exists. 

The conservation laws also validates the accuracy of the proposed algorithms if they keep their initial values even in the nonexistence of the analytical solutions case. The lowest three conservation laws 
\begin{equation}
\begin{aligned}
M&=\int\limits_{-\infty}^{\infty}{udx} \\
E&=\int\limits_{-\infty}^{\infty}{u^2dx} \\
H&=\int\limits_{-\infty}^{\infty}{\dfrac{\mu_1 u^3}{3}+\dfrac{\mu_2 u^4}{6}-\mu_3 (u_x)^2dx}
\end{aligned}
\end{equation} 
are expected to keep their initial values as time proceeds\cite{hamdi1}. In order to measure the absolute relative changes of these quantities at any time $t>0$, $C(M_t)$, $C(E_t)$ and $C(H_t)$ are defined as
\begin{equation}
\begin{aligned}
C(M_t)&=\left | \frac{M_t-M_0}{M_0} \right | \\
C(E_t)&=\left | \frac{E_t-E_0}{E_0} \right | \\
C(H_t)&=\left| \frac{H_t-H_0}{H_0} \right|
\end{aligned}
\end{equation}
where $M_0$, $E_0$ and $H_0$ are initial, $M_t$, $E_t$ and $H_t$ are the quantities at the time $t>0$.

\subsection{Example 1}

\begin{figure}[hp]
	\centering
		\includegraphics[scale =0.75]{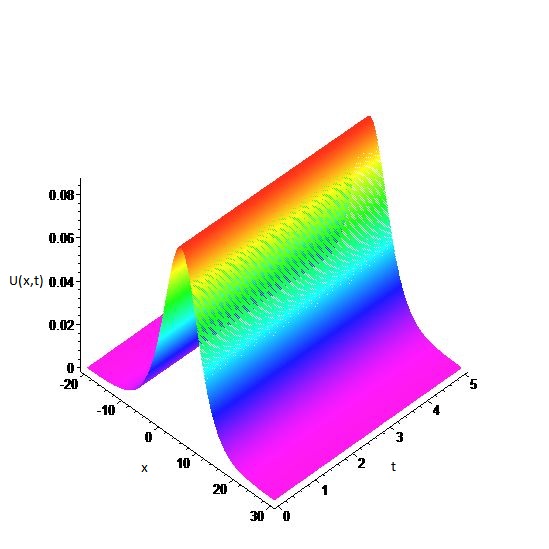}
	\caption{Propagation of initial positive pulse}
	\label{fig:sech}
\end{figure}

The solution representing the propagation of an initial positive pulse is demonstrated by using the initial condition
\begin{equation*}
u(x,0)=\frac{2}{12+3\sqrt{14}\cosh (-\frac{x}{3}+\frac{5}{3})}
\end{equation*}
and the homogeneous Neumann conditions at both ends of the problem interval $[-20,30]$. 
The analytical solution can be written as
\begin{equation*}
u(x,t)=\frac{2}{12+3\sqrt{14}\cosh (-\frac{x}{3}+\frac{5}{3}+\frac{t}{27})}
\end{equation*}
when the compatible parameters  are chosen $\mu _{1}=4,$ $\mu _{2}=-3$ and $\mu _{3}=1$ in the Gardner equation (\ref{gardner}). This solution represents propagation of a positive initial pulse along the $x-$axis as time proceeds. The simulation of the solution is depicted in Fig \ref{fig:sech} in the time domain $[0,5]$ 

\begin{figure}[h]
    \subfigure[Error distribution for $\zeta =1$ at the simulation terminating time]{
   \includegraphics[scale =0.6] {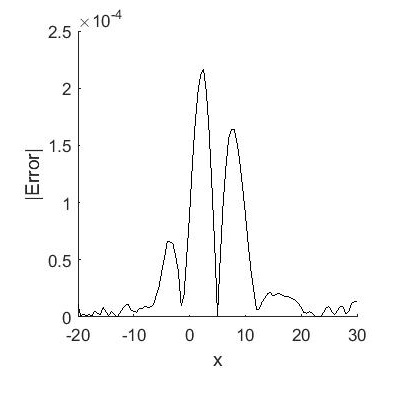}
   \label{fig:1a}
 }
 \subfigure[Error distribution for $\zeta =0.000003$ at the simulation terminating time]{
   \includegraphics[scale =0.6] {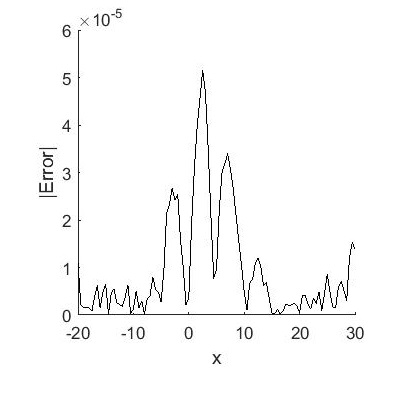}
   \label{fig:1b}
 }
 \caption{Error of the numerical error at the time $t=5$}
\end{figure}

Various discretization parameters $h$ and $\Delta t$ are used to illustrate the numerical solutions. The optimum value of exponential spline parameter $\zeta$ is scanned by comparing the discrete maximum error norms at the simulation terminating time. The maximum values of the errors are measured near the peak points as expected for both $\zeta =1$ and $\zeta =0.000003$. The maximum error distribution for $\zeta =1$ and $\zeta =0.000003$ are depicted in Fig \ref{fig:1a} and Fig \ref{fig:1b}, respectively, for $h=0.5$ and $\Delta t=0.1$. 

The discrete maximum error norms at some distinct times are reported in Table \ref{t1} for various values of the discretization parameter $h$ and fixed $\Delta t$. Even though the solutions are improved when the spatial discretization parameter size is reduced when $\zeta =1$ at both $t=2.5$ and $t=5$. For various values of the spatial discretization parameter, the optimum exponential B-spline parameter is scanned to reduce the maximum error. When $h$ is larger, the optimum $\zeta$ improves the results by decreasing the error approximately one third or one fourth. Even though the better results are obtained for smaller $h$ values, the optimum choice of $\zeta$ increases the accuracy two times.   

\begin{table}[h]{\scriptsize
\caption{$L_{\infty }$ Error norms for $\mu _{1}=4,$ $\mu _{2}=-3,$ $\mu _{3}=1,$ $\Delta t=0.1,$
$t=2.5$ and $5,$ $-20\leq x\leq 30$}
\label{t1}
\begin{tabular}{lllll}
\hline \hline
$N$ & $L_{\infty }(2.5)(\zeta =1)$ & $L_{\infty }(2.5)($various $\zeta $) & $%
L_{\infty }(5)(\zeta =1)$ & $L_{\infty }(5)($various $\zeta$) \\ 
\hline
$100$ & $1.1502\times 10^{-4}$ & $(\zeta =0.000003)3.2331\times 10^{-5}$ & $%
2.1665\times 10^{-4}$ & $(\zeta =0.000003)5.1481\times 10^{-5}$ \\ 
$200$ & $4.1696\times 10^{-5}$ & $(\zeta =0.000001)1.6622\times 10^{-5}$ & $%
5.7428\times 10^{-5}$ & $(\zeta =0.000001)1.8886\times 10^{-5}$ \\ 
$300$ & $2.3860\times 10^{-5}$ & $(\zeta =0.000005)1.3923\times 10^{-5}$ & $%
2.9888\times 10^{-5}$ & $(\zeta =0.000004)1.7006\times 10^{-5}$ \\ 
$400$ & $1.6985\times 10^{-5}$ & $(\zeta =0.000004)1.4470\times 10^{-5}$ & $%
1.8721\times 10^{-5}$ & $(\zeta =0.000003)1.5404\times 10^{-5}$ \\ 
\hline\hline
\end{tabular}}
\end{table}
The initial values of the conservation laws are computed with computer algebra tools and reported in Table \ref{t2}. A stable numerical method is expected to preserve these values as time goes. The absolute relative changes of all related conservation laws reported in Table \ref{t2}. They all are preserved at least six decimal digits for various values of $h$ and $\zeta =1$. The preservation of conservation laws is an indicator of a reliable numerical approach.
\begin{table}[h]{\scriptsize
\caption{Absolute relative changes of conservation laws for for $\zeta =1$}
	\label{t2}
		\begin{tabular}{lllllll}
		\hline \hline
		$N$ & $M_{0}$ & $E_{0}$ & $H_{0}$ & $C(M_{5})$ & $C(E_{5})$ & $C(H_{5}) $ \\ 
\hline
$100$ & $1.04458$ & $0.06013453$ & $0.00407022$ & $5.5668\times 10^{-6}$ & $%
2.6168\times 10^{-8}$ & $1.2174\times 10^{-5}$ \\ 
$200$ & $1.04458$ & $0.06013453$ & $0.00407022$ & $2.9640\times 10^{-6}$ & $%
5.0740\times 10^{-8}$ & $1.0597\times 10^{-6}$ \\ 
$300$ & $1.04458$ & $0.06013453$ & $0.00407022$ & $2.3326\times 10^{-7}$ & $%
2.2152\times 10^{-8}$ & $2.7126\times 10^{-6}$ \\ 
$400$ & $1.04458$ & $0.06013453$ & $0.00407022$ & $1.1862\times 10^{-6}$ & $%
8.8551\times 10^{-10}$ & $3.3555\times 10^{-6}$ \\
	\hline \hline		
		\end{tabular}}
\end{table}

\subsection{Example 2}
\begin{figure}[hp]
	\centering
		\includegraphics[scale =0.75]{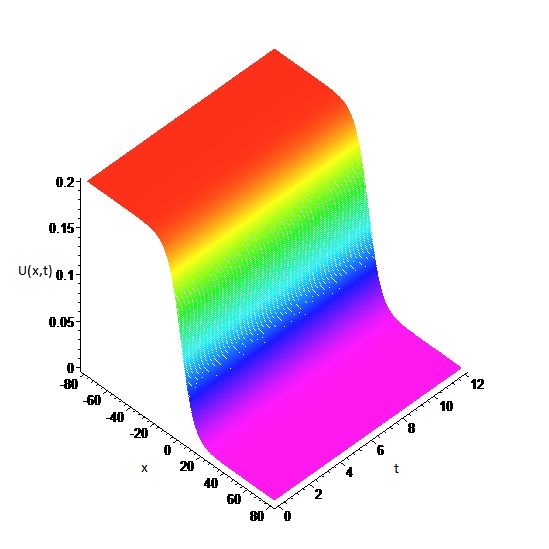}
	\caption{Motion of kink like wave}
	\label{fig:tanh}
\end{figure}

The solution of the Gardner equation (\ref{gardner}) representing the motion of kink like wave is of the form
\begin{equation*}
u(x,t)=\frac{1}{10}-\frac{1}{10}\tanh (\dfrac{(x-\dfrac{t}{30})\sqrt{30}}{60})
\end{equation*}
when the equation parameters are chosen as $\mu_{1}=1, $ $\mu _{2}=-5$ and $\mu _{3}=1$. The initial condition required for the numerical simulation is derived from the analytical solution by substituting $t=0$ in the analytical solution. The homogeneous Neumann boundary conditions are used to complement the problem statement. The designed algorithm is run up to the terminating time $t=12$ in the finite interval $[-80,80]$. The results are graphed in Fig \ref{fig:tanh}. 

\begin{figure}[h]
    \subfigure[Error distribution for $\zeta =1$ at the simulation terminating time]{
   \includegraphics[scale =0.6] {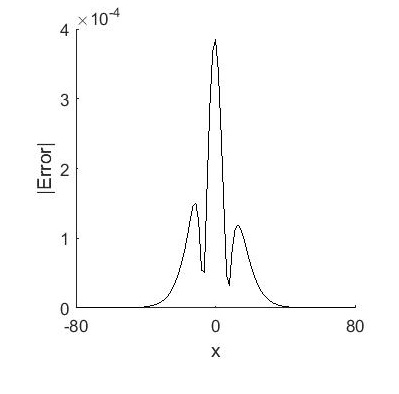}
   \label{fig:2a}
 }
 \subfigure[Error distribution for $\zeta =0.000001$ at the simulation terminating time]{
   \includegraphics[scale =0.6] {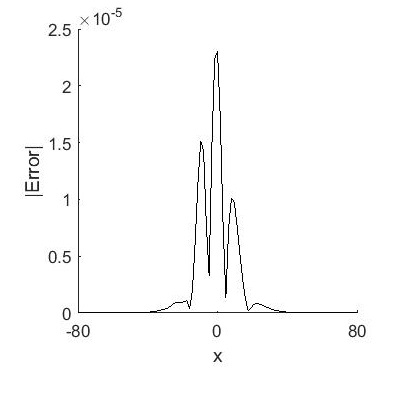}
   \label{fig:2b}
 }
 \caption{Error of the numerical error at the time $t=12$}
\end{figure}

The maximum error distributions for $\zeta =1$ and $\zeta =0.000001$ are depicted in Fig \ref{fig:2a} and Fig \ref{fig:2b}, respectively, for $N=100$. The errors are larger around $x=0$ as expected due to the shape of the kink like wave.

The maximum error norms are calculated for various values of $N$ for $\zeta =1$. In order to improve the results, the optimum value of the exponential B-spline parameter are investigated by scanning. Determination of optimum values of $rho$ different from $1$ gives at least one decimal digit better results, Table \ref{t3}. Even though the results improve when $N$ increases, optimum values of $\zeta$ improves results independent on the discretization parameters.

\begin{table}[h]
\caption{$L_{\infty }$ Error norms for $\mu _{1}=1,$ $\mu _{2}=-5,$ $\mu _{3}=1,$ $\Delta t=0.1,$
$t=12,$ $-80\leq x\leq 80$}
\label{t3}
\begin{tabular}{lllllll}
\hline \hline
$N$ & $L_{\infty}(12)(\zeta =1)$ & $L_{\infty }(12)($various $\zeta$) \\ 
\hline
$100$ & $3.8436\times 10^{-4}$ & $(\zeta=0.000001)2.3022\times 10^{-5}$ \\ 
$200$ &  $1.0016\times 10^{-4}$ & $(\zeta=0.000002)5.8623\times 10^{-6}$ \\ 
$400$ &  $2.5327\times 10^{-5}$ & $(\zeta=0.000004)1.3684\times 10^{-6}$ \\ 
$600$ &  $1.1280\times 10^{-5}$ & $(\zeta=0.000006)5.3420\times 10^{-7}$ \\ 
$800$ &  $6.3476\times 10^{-6}$ & $(\zeta=0.000008)2.3800\times 10^{-7}$ \\
\hline\hline
\end{tabular}
\end{table}
The initial values and absolute relative changes of the conservation laws are tabulated in Table \ref{t4}. Absolute relative changes of all the conservation laws are calculated in four decimal digits independently on number of spatial discretization point number, Table \ref{t4}.
\begin{table}[h]{\scriptsize
\caption{Absolute relative changes of lowest three conservation laws for $\zeta =1$}
	\label{t4}
		\begin{tabular}{lllllll}
		\hline \hline
		$N$ & $M_{0}$ & $E_{0}$ & $H_{0}$ & $C(M_{12})$ & $C(E_{12})$ & $C(H_{12})$ \\ 
\hline
$100$ & $16.1599$ & $3.0129$ & $0.0979$ & $4.9493\times 10^{-4}$ & $%
5.3092\times 10^{-4}$ & $5.4405\times 10^{-4}$ \\ 
$200$ & $16.0799$ & $2.9969$ & $0.0974$ & $4.9750\times 10^{-4}$ & $%
5.3387\times 10^{-4}$ & $5.4720\times 10^{-4}$ \\ 
$400$ & $16.0399$ & $2.9889$ & $0.0971$ & $4.9875\times 10^{-4}$ & $%
5.3531\times 10^{-4}$ & $5.4871\times 10^{-4}$ \\ 
$600$ & $16.0266$ & $2.9862$ & $0.0971$ & $4.9917\times 10^{-4}$ & $%
5.3578\times 10^{-4}$ & $5.4922\times 10^{-4}$ \\ 
$800$ & $16.0199$ & $2.9849$ & $0.0970$ & $4.9937\times 10^{-4}$ & $%
5.3602\times 10^{-4}$ & $5.4947\times 10^{-4}$ \\ 
	\hline \hline		
		\end{tabular}}
\end{table}

\subsection{Example 3}

\begin{figure}[hp]
\subfigure[$t=5$]{
   \includegraphics[scale =0.6] {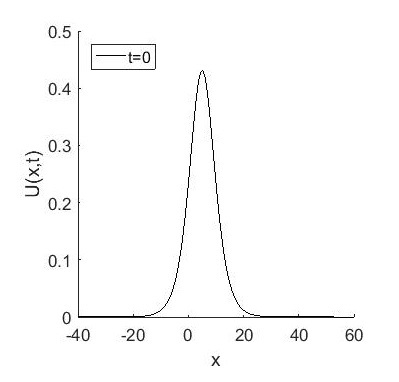}
   \label{fig:3a}
 }
    \subfigure[$t=5$]{
   \includegraphics[scale =0.6] {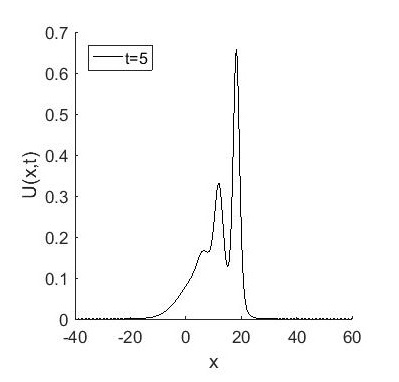}
   \label{fig:3b}
 }
 \subfigure[$t=10$]{
   \includegraphics[scale =0.6] {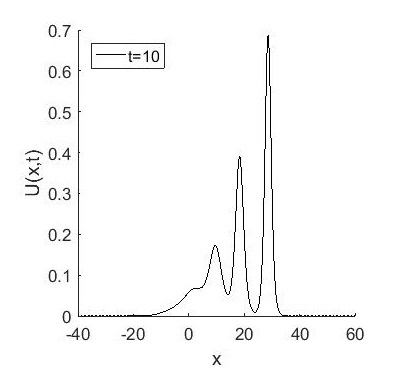}
   \label{fig:3c}
 }
  \subfigure[$t=15$]{
   \includegraphics[scale =0.6] {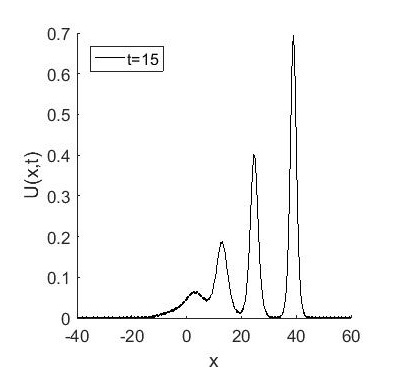}
   \label{fig:3d}
 }
 \caption{Wave generation simulation}
\end{figure}

The wave generation from an initial positive pulse is studied as the last example. The initial condition used in the first example is perturbed carefully to generate new waves. Thus, an initial condition is produced  for $\mu _{1}=10,$ $\mu _{2}=-3$ and $\mu _{3}=1$ as
\begin{equation*}
u(x,0)=\frac{2}{3}\frac{5}{4+\sqrt{14}\cosh (\dfrac{x}{3}-\dfrac{5}{3})}
\end{equation*}
This initial pulse is expected to generate new pulses behind as propagating towards to the right along the horizontal axis. The numerical simulation is accomplished with the discretization parameters $h=0.5$ and $\Delta t=0.1$ in the interval $[-40,60]$. The simulation is depicted in Fig \ref{fig:3a} - \ref{fig:3d}.

The initial values and the absolute relative changes of the conservation laws are recorded during the simulation process duration in Table \ref{t5}. It is clearly observable that all the conservation laws are preserved successfully at least three decimal digits in the simulation period.

\begin{table}[hp]{\scriptsize
\caption{Absolute relative changes of lowest three conservation laws for for $\zeta =1$}
	\label{t5}
		\begin{tabular}{lllllll}
		\hline \hline

$t$ & $M_{0}$ & $E_{0}$ & $H_{0}$ & $C(M_{t})(p=1)$ & $C(E_{t})(p=1)$ & 
$C(H_{t})(p=1)$ \\ 
\hline
$5$ & $5.2255$ & $1.5033$ & $1.5994$ & $1.2040\times 10^{-6}$ & $%
3.7180\times 10^{-5}$ & $21608\times 10^{-3}$ \\ 
$10$ & $5.2255$ & $1.5033$ & $1.5994$ & $3.6819\times 10^{-6}$ & $%
5.2527\times 10^{-5}$ & $3.1907\times 10^{-3}$ \\ 
$15$ & $5.2255$ & $1.5033$ & $1.5994$ & $8.9144\times 10^{-6}$ & $%
5.8526\times 10^{-4}$ & $3.5478\times 10^{-3}$ \\	\hline \hline		
		\end{tabular}}
\end{table}

\section{Conclusion}
The exponential B-spline collocation method is implemented for IBVPs with analytical and non analytical solutions. Having no continuous third order derivative of exponential B-splines forces us to reduce the order of third order derivative term. Thus, the Gardner equation reduces to a coupled nonlinear PDE system with maximum second order derivative term. The exponential B-spline approach to the solution of this system is substituted into the system following the Crank-Nicolson implicit time integration method. Linearization of both non linear terms with quadratic and cubic nonlinearities is followed by implementation of boundary data to obtain a solvable system of equation having equal numbers of equations and unknowns. The algorithm gets ready to run after arranging the initial state by the aid of initial and boundary data.

The validity of the method is checked by computing the error between the analytical and numerical solutions in the case of having an analytical solution. The absolute relative changes of the conservation laws can also be a good indicator to observe the validity and accuracy of the proposed method even when there exits no analytical solution. The scan to determine the optimum exponential B-spline parameter value shows that the accuracy of the solutions can be improved when compared the choice $\zeta =1$.\\
\section{Acknowledgements}
\textsl{This study is supported by Eskisehir Osmangazi University Scientific Research Projects Committee with project number 2016/19052. A brief part was orally presented at 3rd International Conference on Pure and Applied Sciences, Dubai, 2017.}

\end{document}